\documentclass[11pt]{amsart}

\usepackage[utf8]{inputenc}
\usepackage[T1]{fontenc}

\usepackage{tcolorbox}
\usepackage{amssymb}
\usepackage{amsmath}
\usepackage{amsthm}

\usepackage{graphicx}

\usepackage{stmaryrd}
\usepackage{mathrsfs}
\usepackage{amscd}
\usepackage{bbm}
\usepackage{enumitem}
\usepackage{fncylab}
\usepackage{xcolor}
\usepackage{url}
\usepackage{leftidx}
\usepackage{calc}
\usepackage[width=\textwidth]{caption}
\usepackage[normalem]{ulem}
\usepackage{dsfont}

\definecolor{boubelcolor}{rgb}{.65,0.05,0}

\DeclareUnicodeCharacter{2260}{\colorlet{memoireboubel}{.}\color{boubelcolor}} 
\DeclareUnicodeCharacter{00B1}{\color{memoireboubel}{}} 
\DeclareUnicodeCharacter{00A3}{\colorlet{memoirejuillet}{.}\color{blue}} 
\DeclareUnicodeCharacter{00A7}{\color{memoirejuillet}{}} 

\newtheorem*{remark}{Remark}

\theoremstyle{plain}

\theoremstyle{definition}

\newtheorem{notation-terminology}[them]{Notation/Terminology}

\newtheorem{reminder-notation}[them]{Reminder/Notation}

\theoremstyle{remark}

\newcommand{\dd}{\mathrm{d}}

\renewcommand{\leq}{\leqslant} 
\renewcommand{\geq}{\geqslant} 

\newtheorem{theorem}{Theorem}

\newcommand\longoverline{\bgroup\markoverwith
{{\rule[1.5ex]{2pt}{0.4pt}}}\ULon}  

\newcounter{point}
\newcounter{souspoint}[point]

\labelformat{point}{(#1)}          
\labelformat{souspoint}{(#1)}  

\renewcommand{\P}{\mathbb{P}}
\newcommand{\p}{\mathcal{P}}

\newcommand{\F}{\mathcal{F}}
\newcommand{\G}{\mathcal{G}}

\newcommand{\R}{\mathbb{R}}

\newcommand{\E}{\mathbb{E}}

\newcommand{\cx}{\prec_{\mbox{\tiny cx}}}
\newcommand{\icx}{\prec_{\mbox{\tiny icx}}}
\newcommand{\st}{\prec_{\mbox{\tiny st}}}

\title{A coupling proof of convex ordering for compound distributions}
\author{Jean B\'erard and Nicolas Juillet}

\begin{document}
\begin{abstract}

In this paper, we give an alternative proof of the fact that, when compounding a nonnegative probability distribution, convex ordering between the distributions of the number of summands implies convex ordering between the resulting compound distributions. Although this is a classical textbook result in risk theory, our proof exhibits a concrete coupling between the compound distributions being compared, using the representation of one-period discrete martingale laws as a mixture of the corresponding extremal measures.

\end{abstract}
\maketitle

\section{Introduction}

 Consider an i.i.d. sequence of nonnegative random variables $\mathscr{X}=(X_i)_{i \geq 1}$, and two integer-valued random variables $M,N$, independent from the sequence $\mathscr{X}$. Assume that $\E(X_1)<+\infty$, $\E(M)<+\infty$,  $\E(N)<+\infty$, and that a comparison between $M$ and $N$ holds with respect to the convex\footnote{We refer to  \cite{ShaSha, MulSto} for the definition and main properties of the convex ordering of probability measures.} ordering:  $M \cx N$. We then have the following comparison between the compound variables $X_1+\cdots + X_M$ and $X_1+\cdots+X_N$ with respect to the convex ordering:
\begin{equation}\label{e:comp-convexe}X_1+\cdots + X_M \cx X_1+\cdots+X_N.\end{equation}
This is a classical result (see e.g. Theorem 4.A.9 in \cite{ShaSha} or Theorem 4.3.6 in  \cite{MulSto})\footnote{In \cite{ShaSha,MulSto}, the theorems are stated with respect to the increasing convex order, and involve two sequences of random variables instead of just one, so that \eqref{e:comp-convexe} appears as a special case of these results.  We refer to Section \ref{ss:other-orders} for a discussion of how the more general case can be deduced from \eqref{e:comp-convexe}.}, useful in the context of risk theory (see e.g. \cite{KaaDen}, chap. 7) where its interpretation is that compounding with a riskier frequency distribution leads to a riskier aggregated loss distribution.

The proof given in \cite{ShaSha} is analytical in nature, and consists in showing that, given a non-decreasing convex function $f  \ :  \ [0,+\infty[ \to  \R$, the sequence $(u_n)_{n \geq 0}$ defined by $u_n = \E f(X_1+\cdots+X_{n})$ satisfies, for all $n \geq 0$, the condition $u_{n+2}-u_{n+1} \geq u_{n+1}-u_n$ (provided that $u_{n+2}, u_{n+1}, u_n$ have a finite value)\footnote{The proof given in \cite{MulSto} is similar, using functions of the special form $f(x)=(x-t)^{+}$.}. Here, we give an alternative proof of \eqref{e:comp-convexe}, based on a coupling between the two random variables being compared, which provides an explicit realization of Strassen's condition (see \cite{Str} or e.g. \cite{MulSto} section 1.5.2) for convex ordering: we construct a pair of random variables $(A,B)$ such that 
\begin{equation}\label{e:Strassen-1} A \stackrel{d}{=} X_1+\cdots + X_M, \ B \stackrel{d}{=} X_1+\cdots + X_N, \mbox{ and } A = \E(B | \sigma(A)) \mbox{ a.s.}\end{equation}

Our proof relies a representation result which we call a {\it diatomic representation} of convex ordering, stated as Theorem \ref{t:diatomique} in Section \ref{s:diatomic-representation}, where we review several approaches for proving the existence of this representation, including an explicit algorithm in the case of discrete distributions. Section \ref{s:coupling-construction} contains the coupling construction leading to \eqref{e:Strassen-1} and the proof of \eqref{e:comp-convexe}. Finally,  in Section \ref{s:final-remarks-extensions}, we discuss various extensions of these results.

\section{Diatomic representation}\label{s:diatomic-representation}

Given two real numbers $x<y$, we define, for all $z$ the normalized barycentric coordinates:
\[\alpha^{x,y}(z) = \frac{y-z}{y-x} \mbox{ and }   \ \beta^{x,y}(z) = 1- \alpha^{x,y}(z) = \frac{z-x}{y-x}.\]
In the case $x=y$, we extend the above definition by setting $\alpha^{x,y}(z)=1$ and $\beta^{x,y}(z)=0$. Whenever $x \leq z \leq y$, both $\alpha^{x,y}(z)$ and $\beta^{x,y}(z)$ lie in the interval $[0,1]$, and $z$ can be written as the convex combination of $x$ and $y$ : 
$$z=\alpha^{x,y}(z)\cdot x+\beta^{x,y}(z)\cdot y.$$

\begin{theorem}\label{t:diatomique}
Given two probability distributions $\mu,\nu$ on $\R$ possessing a finite expectation, the comparison $\mu \cx \nu$ holds if and only if there exists a triple of random variables $(V_-,U,V_+)$ defined on the same probability space and such that:
\begin{equation}\label{e:dans-l-ordre}U \in [ V_-, V_+]\end{equation}
\begin{equation}\label{e:marginale-U} \mbox{Law}(U)= \mu\end{equation}
\begin{equation}\label{e:di-marginale} \nu=\E \left[ \alpha^{V_-, V_+}(U) \cdot \delta_{V_-} +  \beta^{V_-, V_+}(U) \cdot \delta_{V_+}  \vphantom{\sum}  \right] \end{equation}
We call such a triple $(V_-,U,V_+)$ a {\it diatomic representation} of the stochastic ordering  $\mu \cx \nu$.
\end{theorem}
A more concrete statement of \eqref{e:di-marginale} is that   $\mbox{Law}(V) = \nu$, where $V$ is a random variable whose conditional distribution  with respect to $V_-,U,V_+$ is given by:
\begin{align*}
V=
\begin{cases}
V_-&\text{with probability }\alpha^{V_-, V_+}(U)\\
V_+&\text{with probability }\beta^{V_-, V_+}(U)
\end{cases}
\end{align*}

Theorem \ref{t:diatomique} may not have been stated under this specific form in the mathematical literature, but its content is certainly not new. In the following subsections we review several ways of proving this result.

\subsection{Proof of Theorem \ref{t:diatomique} via Choquet's and Douglas' theorems (for compactly supported measures $\nu$)}

We prove here the result only in the case of measures $\mu\cx \nu$ concentrated on some closed interval $K$. 

Consider the following space of measures on $\R^2$: 
\[S_K=
\left\{\pi\in\mathcal{M}_+(K^2):\,
\forall a\in \R,\,\forall b\in \mathcal{C}_b(\R),\,\iint a+b(u).(v-u)\dd\pi(u,v)=a
\right\}.
\]
Here $\mathcal{M}_+(K^2)$ is the space of finite positive Borel measures on $\R^2$ with support contained in $K^2$, and $\mathcal{C}_b(\R)$ the space of real-valued continuous and bounded functions on $\R$. In other words, denoting by $F$ the space of functions $f_{a,b}:(u,v)\in \R^2\to a+b(u)(v-u)$, the above definition reads:  $$S_K=\{\pi\in\mathcal{M}_+(K^2):\,\forall f_{a,b}\in F,\,\iint f_{a,b}\dd \pi=a \}.$$ It turns out that $S_K$ is a non-empty space of probability measures, known as martingale measures in the literature since for $(X_1,X_2)\sim \pi\in S_K$, the process $(X_i)_{i=1,2}$ is a martingale on its natural filtration.

Our task amounts to proving that a measure $\pi\in S_K$ is represented as a mixture of `triplet' measures of the form $\delta_{z}\otimes [\alpha \delta_x+\beta \delta_y]$, where $z=\alpha x+\beta y$ and $\alpha,\,\beta$ are nonnegative and satisfy $\alpha+\beta=1$.

We admit without proof that $S_K$ is convex and compact for the weak topology. Choquet's Theorem (see \cite{Ch} or e.g. \cite{Phe}) then implies that every measure in $S_K$ can be represented as a mixture of extremal measures in $S_K$. So we shall be done as soon as we can prove that the extremal measures in $S_K$ are triplet measures concentrated on $S_K$. Let $\eta$ be such an extremal measure and $\mu,\,\nu$ its marginals. We aim at proving that $\mu$ is supported on a single point and that $\nu$ is supported on at most two points. Striving for a contradiction suppose that  $\mu$ is not a Dirac measure. Hence, there exists a Borel set $E$ such that $\mu(E)\notin \{0,1\}$. We consider the $L^1(\eta)$ distance between $f:(u,v)\to\mathds{1}_E(u)$ and the functions $f_{a,b}\in F$. Letting $(U,V)$ be distributed according to $\eta$ we find
\begin{align*}
\|f-f_{a,b}\|_{L^1(\eta)}&\geq \E\left(\left| \E(f(U,V)-f_{a,b}(U,V)|\,U)  \vphantom{\sum} \right| \right)\\
&\geq \E \left( \left| \E(\mathds{1}_E(U)-a|\,U)   \vphantom{\sum}  \right| \right)\\
&= \P(U\in E)|1-a|+\P(U\notin E)|a|\\
&\geq 1/2 \min(\P(U\in E),\P(U\notin E)).
\end{align*}
Note that this lower bounds only depends on $f$. We have thus proved that $F$ is not dense in $L^1(\eta)$. According to Douglas's theorem \cite{Dou} (see also \cite[Chapter V]{ReYo99}) this in contradiction with the fact that $\eta$ is extremal. Thus $\eta$ is of type $\delta_{u_0}\otimes \nu$ where $u_0$ is the barycenter of $\nu$, i.e $\int v\,\dd\nu(v)=u_0$. Next, again striving for a contradiction suppose that there exists a partition $(A_i)_{i\in\{1,2,3\}}$ of $\R$ such that $\nu(A_i)>0$ holds for every $i$. From Douglas's Theorem again the set $F$ of functions is dense in $L^1(\eta)$. In particular any function $g_{c_1,c_2,c_3} : (u,v)\mapsto \sum_{i=1}^3c_i\mathds{1}_{A_i}(v)$ can be approximated in $L^1(\eta)$ by functions $f_{a,b}$. The linear map $(c_1,c_2,c_3)\mapsto (\int_{A_1} g_{c_1,c_2,c_3}\,\dd\eta, \int_{A_2} g_{c_1,c_2,c_3}\,\dd\eta, \int_{A_3} g_{c_1,c_2,c_3}\,\dd\eta)$ is clearly linear and onto. It follows that the linear map 
\[(a,b)\in \R^2\mapsto (\int_{A_1} f_{a,b}\,\dd\eta, \int_{A_2} f_{a,b}\,\dd\eta, \int_{A_3} f_{a,b}\,\dd\eta)\in \R^3\]
 is onto as well, a contradiction. Therefore, extremal measures of $S_K$ are of type $\delta_{z}\otimes [\alpha \delta_x+\beta \delta_y]$ where $z=\alpha x+\beta y$, as it was required.

\subsection{Proof of Theorem \ref{t:diatomique} via Strassen's theorem (general case)}\label{s:second_proof}
Another approach to proving Theorem \ref{t:diatomique} is to use Strassen's theorem instead of Choquet's and Douglas's theorems: there exists a kernel $k:\R\to \p(\R)$ such that $\mu$ almost surely $k^u=k(u)$ has mean $u$ and it holds $\mu\cdot k=\nu$. (Such kernels are known as dilations or martingale kernels in the literature.) Hence the mixture with weight $\mu$ of the measures $\delta_u\otimes k^u$ defines a probability measure $\pi$ on $\R^2$ whose marginals are $\mu$ and $\nu$. To complete the proof, it remains to check that each measure $k^u$ can be represented as a mixture of diatomic measures with mean $u$. This last fact is a classical step in the proof of Skorokhod's representation theorem (see e.g.~\cite[Theorem 8.1.1]{Dur}): every probability measure on $\R$ with mean $u$ can be represented as a mixture of diatomic measures with mean $u$. See \cite[\S 5.1]{Ju_seminaire} for another approach.

\begin{remark}
The search for martingale kernels $k:u\mapsto k^u$ is a key question in the field of martingale optimal transport.
The first completely canonical method seems to be the left-curtain coupling by Beiglböck and Juillet \cite{BJ} that is also of particular interest to us. Not only is $k^u$  canonical but when $\mu$ is diffuse its kernels $k^u$ are automatically diatomic (this also holds for the former coupling by Hobson and Neuberger \cite{HoNe12} under more general assumptions). This is not the case  if $\mu$ possesses atoms. However a quantile version of the left-curtain coupling is described in a second paper by the same authors \cite{BJ2} where the martingale measure $\pi$ directly appears as a mixture over the set $[0,1]$ of quantile levels $\omega\in [0,1]$ of diatomic kernels $\delta_{z_\omega}\otimes(\alpha^{x_\omega,y_\omega}(z_\omega)\delta_{x_\omega}+\beta^{x_\omega,y_\omega}(z_\omega)\delta_{y_\omega})$ where $z_\omega$ is the $\omega$-quantile of $\mu$. Note that the same can be said of the recent coupling by Jourdain and Margheriti \cite{JM}.
\end{remark}

\subsection{Algorithmic proof of Theorem \ref{t:diatomique} (for finitely supported $\mu,\,\nu$).}\label{ss:algo}

\begin{figure}

\center{\includegraphics[width=12cm]{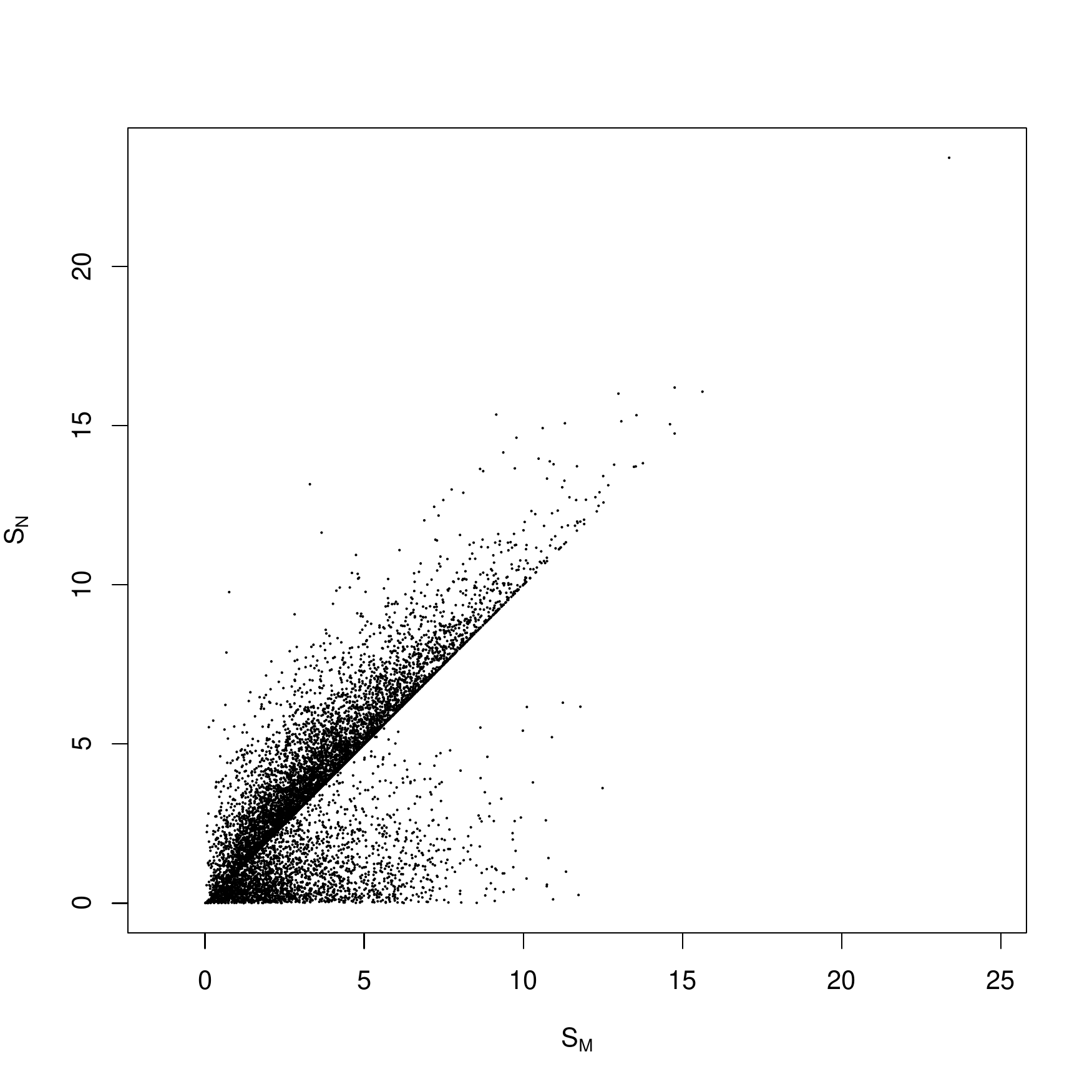}}
\caption{The figure shows $10^4$ simulated pairs $(A,B)$, with $\mu= \sum_{n=2}^5 \textstyle{\frac{1}{4}} \delta_n$ and $\nu= \textstyle{\frac{1}{4}} \delta_1 +  \textstyle{\frac{1}{4}} \delta_3 + \sum_{n=4}^6 \textstyle{\frac{1}{6}} \delta_n$, $\mbox{Law}(X_i) = \mathscr{E}(1)$, and the algorithm of Section \ref{ss:algo} to produce the diatomic decomposition.}
\label{f:image-couplage}
\end{figure}

We now describe an explicit algorithmic construction, inspired by \cite{BJ,BJ2}, leading to a diatomic decomposition in the case where both $\mu$ and $\nu$ are finitely supported probability measures. This algorithm is used to produce the simulation shown in Fig. \ref{f:image-couplage}.

Let us write $\mu=\sum_{i=1}^p \mu(u_i) \delta_{u_i}$ and  $\nu=\sum_{j=1}^q \nu(v_j) \delta_{v_j}$. 

\begin{tcolorbox}
{\bf Initialization:}
$\mu_{*} \leftarrow \mu $, $\nu_{*} \leftarrow \nu $, $\theta_{*}   \leftarrow 1 $, $\mathscr{T} \leftarrow \emptyset$
 \end{tcolorbox}
\begin{tcolorbox}
{\bf Loop:} Repeat the following steps while $\theta_{*} > 0$ 

Pick a triple $(v_{j-},u_i,v_{j+})$ such that

\begin{itemize}
\item[a.] $\nu_{*}(v_{j-})>0$, $\mu_{*}(u_i)>0$, $\nu_{*}(v_{j+})>0$

\item[b.] $v_{j-} \leq u_i \leq v_{j+}$

\item[c.] $\nu_{*}( (v_j-,v_j+) )=0$

\end{itemize}
$s \leftarrow  \min\left(  \mu_{*}(u_i), \nu_{*}(v_{j-})/ \alpha^{v_{j-},v_{j+}}(u_i) , \nu_{*}(v_{j+})/ \beta^{v_{j-},v_{j+}}(u_i))        \right)$

 $\mu_{*} \leftarrow \mu_{*} - s \delta_{u_i}$

$\nu_{*} \leftarrow \nu_{*} - s \alpha^{v_{j-},v_{j+}}(u_i) \delta_{v_{j-}} - s \beta^{v_{j-},v_{j+}}(u_i) \delta_{v_{j+}}$

 $\theta_{*} \leftarrow \theta_{*} -s $

$\mathscr{T} \leftarrow \mathscr{T} \cup \{ \left( (v_{j-} , u_i , v_{j+}),  s  \right)  \}$

\end{tcolorbox}

\begin{tcolorbox}
{\bf Result:} return the set $\mathscr{T}$
\end{tcolorbox}

The probability distribution of $(V-,U,V+)$ is then deduced from $\mathscr{T}$ as $$\sum_{(v-,u,v+,s) \in \mathscr{T} } s \delta_{(v-,u,v+)}.$$

The reason why the above algorithms stops lies in the fact that the transformation performed on $\mu_{*}$ and $\nu_{*}$ keeps the comparison $\mu_{*} \cx \nu_{*}$ valid throughout the execution of the algorithm
(see the proof of Lemma 2.8 in \cite{BJ}), with $\mu_{*}$ and $\nu_{*}$ having equal total mass.
As a consequence, as long as $\mu_{*}$ and $\nu_{*}$ do not have zero total mass, the comparaison $\mu_{*} \cx \nu_{*}$  ensures that a triple $(v_{j-},u_i,v_{j+})$ satisfying conditions a.-b.-c. exists. Finally, since at each step at least one of the three numbers $\nu_{*}(v_{j-}),\mu_{*}(u_i),\nu_{*}(v_{j+})$ is set to zero, the total mass of both $\mu_{*}$ and $\nu_{*}$ must reach zero after a finite number of steps.

\begin{remark}
If, in the loop part of the algorithm, one systematically choses the unique triple $(v_-,u,v_+)$ such that $u$ is the leftmost point in the support of $\mu_{*}$ (such a choice is always possible, see the proof of Lemma 2.8 in \cite{BJ}), the end-result of the algorithm is the so-called left-curtain coupling. 
\end{remark}

\section{Coupling construction}\label{s:coupling-construction}

We now describe the coupling construction leading to our proof of \eqref{e:comp-convexe}. Consider a triple $(N_-,M,N_+)$ as in Theorem \ref{t:diatomique}, with $\mu=\mbox{Law}(M)$ and $\nu=\mbox{Law}(N)$, and an i.i.d. sequence $\mathscr{X}=(X_i)_{i \geq 1}$ independent from $(N_-,M,N_+)$. For all integer $k \geq 1$, we let $S_k=\sum_{i=1}^k X_i$, with the convention that $S_0 = 0$.  Finally, we let $\F= \sigma(M,N_-,N_+, S_{N_-},S_{N_+})$ and $\G= \sigma(M,N_-,N_+, S_{M},S_{N_-},S_{N_+})$. Note that $\G=\F\vee \sigma(S_M)$.

We start our construction by setting:
\begin{equation}\label{e:def-A}A =S_M.\end{equation}
Next, we specify $B$ by the requirement that, conditional upon $\G$, the distribution of $B$ is: 
\begin{equation}\label{e:def-B}\alpha^{S_{N_-},S_{N_+}}(A) \cdot \delta_{ S_{N_-} }  + \beta^{S_{N_-},S_{N_+}}(A)  \cdot \delta_{ S_{N_+} }.\end{equation}
Note that $\alpha^{S_{N_-},S_{N_+}}(A)$ and $\beta^{S_{N_-},S_{N_+}}(A)$ do indeed lie in the interval $[0,1]$ thanks to the assumption that the random variables $X_i$ are nonnegative, so that $S_{N_-} \leq A = S_M \leq S_{N_+}$, since $N_- \leq M \leq N_+$.

We now proceed to checking that all three properties listed in \eqref{e:Strassen-1} are satisfied by the random variables $A$ and $B$. From \eqref{e:def-A}, it is immediate that $A$ has the required distribution. Moreover, from the definition of $\alpha$ and $\beta$, one has the identity 
$$\alpha^{S_{N_-},S_{N_+}}(A) \cdot S_{N_-} +   \beta^{S_{N_-},S_{N_+}}(A) \cdot  S_{N_+}  = A,$$
which rewrites as:
$$\E \left[  B |  \G  \right] = A \mbox{ a.s.},$$
whence, taking the conditional expectation $\E(\cdot | \sigma(A))$ on both sides, and using the fact that $\sigma(A) \subset \G$ since $A=S_M$ is $\G-$measurable, 
$$\E \left[  B |  \sigma(A)  \right] = A \mbox{ a.s.},$$
as required by \eqref{e:Strassen-1}.

To conclude the proof, it remains to check that $B \stackrel{d}{=} S_N$. By construction, the conditional distribution of $B$ given $\F$ can be written as: 
$$ \E \left[\alpha^{S_{N_-},S_{N_+}}(A) \left| \vphantom{\sum} \right.  \F \right] \cdot  \delta_{ S_{N_-} }  + \E \left[ \beta^{S_{N_-},S_{N_+}}(A) \left| \vphantom{\sum} \right.    \F   \right] \cdot  \delta_{ S_{N_+} }.$$
Now observe that, by symmetry, given integers $n_- \leq m \leq n_+$ such that $n_- < n_+$, we have: 
$$\E\left[ S_m - S_{n_-} \left| \sigma(S_{n_-},S_{n_+} ) \vphantom{\sum} \right. \right] = \frac{m - n_-}{n_+-n_-} (S_{n_+}-S_{n_-}) \mbox{ a.s.},$$
from which we deduce that\footnote{In the case where $n_-=m=n_+$, these identities are still (obviously) valid.}
$$ \E \left[\alpha^{S_{n_-},S_{n_+}}(S_m)  \left| \vphantom{\sum} \right. \sigma(S_{n_-},S_{n_+} ) \right]=   \alpha^{n_-,n_+}(m) \mbox{ a.s.}$$
and 
$$ \E \left[\beta^{S_{n_-},S_{n_+}}(S_m)  \left| \vphantom{\sum} \right. \sigma(S_{n_-},S_{n_+} )\right]=   \beta^{n_-,n_+}(m) \mbox{ a.s.}$$
As a consequence, the conditional distribution of $B$ given $\F$ is none but: 
$$ \alpha^{N_-,N_+}(M) \cdot \delta_{S_{N_-}} +  \beta^{N_-,N_+}(M)  \cdot \delta_{S_{N_+}}.$$
From the fact that the sequence $\mathscr{X}=(X_i)_{i \geq 1}$ is independent from $(N_-,M,N_+)$, and from condition \eqref{e:di-marginale}, we deduce that 
$$B \stackrel{d}{=} S_N,$$
which concludes the proof.
\section{Final remarks and extensions}\label{s:final-remarks-extensions}

\subsection{Continuous time}\label{ss:cont}
We note that the coupling construction described in Section \ref{s:coupling-construction} can be extended in continuous time. For instance, let $\mathscr{N}=(N_t)_{t\geq 0}$ be a standard Poisson process and $S\cx T$ nonnegative integrable random variables independent from $\mathscr{N}$. Then one has $N_S\cx N_T$, and it is straightforward to extend our approach to define a pair of random variables $(A,B)$ such that \begin{equation}\label{e:Strassen} A \stackrel{d}{=} N_S, \ B \stackrel{d}{=} N_T, \mbox{ and } A = \E(B | \sigma(A)) \mbox{ a.s.}\end{equation}
The same approach still works in exactly the same way if we consider an integrable subordinator instead of a Poisson process.

\subsection{Exchangeable random variables}\label{ss:exch}
If the sequence of random variables $(X_i)_{i \geq 1}$ is assumed to be exchangeable instead of i.i.d., the coupling construction described in Section \ref{s:coupling-construction} works in exactly the same way. (The classical proof found in \cite{ShaSha, MulSto} also works in this case.) 
Note that, in the case of an infinite exchangeable sequence of random variables, one can directly deduce \eqref{e:comp-convexe} from the i.i.d. case, using the De Finetti representation of such a sequence as a mixture of (distributions of)  i.i.d. sequences, and the characterization of \eqref{e:comp-convexe} through the inequality \begin{equation}\label{e:fonc-convexe}\E f(X_1+\cdots+X_M) \leq \E f(X_1+\cdots+X_N)\end{equation} for all convex functions $f$. On the other hand, if $M$ and $N$ are assumed to have finite support, say $\{0,1,\ldots, q\}$, and one considers a finite exchangeable sequence of random variables $X_1,\ldots, X_q$, the extension of De Finetti's theorem to this case (see \cite{KerSze,JanKonYua}) leads in general to a signed mixture of i.i.d. sequences, so one cannot integrate the inequality \eqref{e:fonc-convexe} with respect to the mixing measure in order to directly deduce \eqref{e:comp-convexe}.  

\subsection{Increasing convex ordering}\label{ss:other-orders}
Assume that a comparison between $M$ and $N$ holds with respect to the increasing convex ordering: $M \icx N$. We then have the following modified version of \eqref{e:comp-convexe}:  \begin{equation}\label{e:comp-convexe-croissant}X_1+\cdots + X_M \icx X_1+\cdots+X_N.\end{equation}

To deduce \eqref{e:comp-convexe-croissant} from \eqref{e:comp-convexe}, we note that the comparison $M \icx N$ implies that there exists an integer-valued\footnote{The existence of a random variable $N_0$ such that $M \st N_0$ and $N_0 \cx N$ is a classical and easily proved decomposition result for the increasing convex order. That $N_0$ can, in addition, be chosen to be integer-valued is less standard. One possible proof of this fact is that, when $\mu \icx \nu$, there exists a kernel $k:\R\to \p(\R)$ such that $\mu$ almost surely $k^u=k(u)$ has mean $ \geq u$ and it holds $\mu\cdot k=\nu$. In turn, $k^u$ can be written as a mixture $k^u = \alpha_u k^u_1 + (1-\alpha_u)k^u_2$, where $\mu$ almost surely $\alpha_u \in [0,1]$, $k^u_1$ has mean $u$, and the support of $k^u_2$ is contained in $[u,+\infty[$. For $k'(u)=\alpha_u  \delta_u+ (1-\alpha_u)k^u_2$ we have $\mu\st (\mu\cdot k')\cx \nu$ and a coupling of the corresponding random variables can also be easily deduced. Alternatively for another proof one can consider Kellerer's kernels defined in \cite[\S 2.1]{Ke72} for connecting  $\mu \icx \nu$.} random variable $N_0$ such that $M \st N_0$ and $N_0 \cx N$, where $\st$ denotes the usual stochastic ordering. Given an i.i.d. sequence  $\mathscr{X}=(X_i)_{i \geq 1}$ independent from $M,N_0,N$, the fact that the $X_i$s are nonnegative random variables, combined with $M \st N_0$, yields the comparison $X_1+\cdots+X_M \st  X_1+\cdots+X_{N_0}$. Then, using  \eqref{e:comp-convexe}, we deduce that $X_1+\cdots+X_{N_0} \cx X_1+\cdots+X_{N}$, so that \eqref{e:comp-convexe-croissant} is proved.

Now consider two sequences $\mathscr{X}=(X_i)_{i \geq 1}$ and $\mathscr{Y}=(Y_i)_{i \geq 1}$ of i.i.d. nonnegative random variables, and, in addition to $M \icx N$, assume that $X_i \icx Y_i$. We then have the following extension of \eqref{e:comp-convexe-croissant}:
\begin{equation}\label{e:comp-convexe-croissant-double}X_1+\cdots + X_M \icx Y_1+\cdots+Y_N,\end{equation}
which is a straightforward consequence of \eqref{e:comp-convexe-croissant} and of the fact that the increasing convex ordering is preserved by convolution.

\subsection{Counterexample to \eqref{e:comp-convexe} when the $X_i$ are not positive random variables}  
We let $M \equiv 1$ and $N \sim \frac12(\delta_{0}+\delta_2)$, so that the condition $M \cx N$ is satisfied. Then let $X_1,X_2$ be i.i.d. random variables, independent from $N$ (and $M$), with 
$X_i \sim \frac12(\delta_{-1}+\delta_1)$. We find that $X_M \sim \frac12(\delta_{-1}+\delta_1)$, while $X_N \sim \frac18(\delta_{-2}+6\delta_0+\delta_2)$. Hence $1=\E(|X_N|)<\E(|X_M|)=\frac12$ so that $X_M\cx X_N$ cannot hold.

\bibliographystyle{plain}

\bibliography{comparaison-convexe}

\end{document}